\documentclass[english,12pt,a4paper]{amsart}
\usepackage[T1]{fontenc}
\usepackage{lmodern}
\usepackage[width=14cm, height=21cm]{geometry}
\usepackage{graphicx}
\usepackage{mathtools}
\usepackage{amssymb}
\usepackage{amsthm}
\usepackage{babel}
\usepackage{hyperref}
\usepackage{enumerate}
\usepackage{amsmath}

\theoremstyle{plain}
\newtheorem{theorem}{Theorem}[section]

\newtheorem{lemma}{Lemma}[section]

\begin{document}
	\title[Ramanujan's Continued Fractions of Orders Five, Ten, and Twenty]{On Ramanujan's Continued Fractions of Orders Five, Ten, and Twenty and Associated  Lambert Series Identities}
\author{Shruthi C. Bhat}
\address{Shruthi C. Bhat, Manipal Institute of Technology, Manipal Academy of Higher Education, Manipal-576104, Karnataka, India.}
\email{shruthi.research1@gmail.com}

\author{B. R. Srivatsa Kumar}
\address{B. R. Srivatsa Kumar, Manipal Institute of Technology, Manipal Academy of Higher Education, Manipal-576104, Karnataka, India.}
\email{srivatsa.kumar@manipal.edu}
	\maketitle
	\begin{abstract}
		 In this work, we establish several  new identities connecting Ramanujan's continued fractions of order twenty. By employing product representation for Jacobi's theta function $\theta_1$, we derive a family of new relations connecting the continued fractions of order twenty with continued fractions of order ten and Rogers-Ramanujan continued fraction. Further, utilizing certain mock theta functions and their logarithmic derivatives, we obtain beautiful relations between Lambert series and theta functions of level twenty. Using Ramanujan's  $_1 \psi_1$ summation formula, we establish Lambert series identities associated with the continued fractions of order twenty. These results extend earlier work on continued fractions of order 6, 12, and 16 and contribute to   theory of $q$-series.
	\end{abstract}
	

	\section{Introduction}
The theory of continued fractions introduced by Ramanujan occupies key role in the study of $q$-series, theta functions, and modular forms. It is interesting to study the modular properties of his continued fractions through the product expansions. A detailed description of Ramanujan's continued fractions and $q$-series identities can be found in his notebooks \cite{Ramanujan1957, Ramanujan2000} and in extensive commentaries of Berndt \cite{Berndt1991}. In the recent years, many mathematicians found interest in the theory of continued fractions and revisited many of his continued fractions.   
The following is the renowned Rogers-Ramanujan continued fraction \cite{Ramanujan2000} 
\begin{align}\label{R}
	R(q):=q^{1/5} \dfrac{f(-q,-q^4)}{f(-q^2,-q^3)}=\dfrac{q^{1/5}}{1+\dfrac{q}{1+\dfrac{q^2}{1+\ldots}}},
\end{align}
where	
\begin{align}\label{id1}
	f(x,y)=\sum_{n=-\infty}^{\infty} x^{n(n+1)/2}y^{n(n-1)/2},  \quad |xy|<1,
\end{align} is the renowned Ramanujan's general theta function \cite[p. 35, Entry 19]{Berndt1991}. $f(x,y)$ is one of the examples of $q$-series that has both infinite sum and infinite product representations. 

For the wonderful work on $R(q)$, one can refer  \cite{Ramanujan2000, Berndt1991, Andrews2005}. Inspired by Ramanujan's work, Liu \cite{Liu2007} and Chan \textit{et al.} \cite{Chan_2009} proved many new identities associated to	  $R(q)$. Also, they established new Eisenstein series identities involving $R(q)$. Further, Cao \textit{et al.} \cite{Cao_2023} obtained a  theta function identity which resembles Rogers-Ramanujan-Slater type identity.
The following continued fractions of order ten were studied by Rajkhowa and Saikia \cite{Rajkhowa_2025}, where algebraic relations and explicit evaluations of these continued fractions were obtained.
\begin{align}\label{S1}
	\nonumber	S_1(q):&=q^{3/4} \dfrac{f(-q,-q^{9})}{f(-q^4,-q^6)} \\
	&=\dfrac{q^{3/4}(1-q)}{(1-q^{5/2})+\dfrac{q^{5/2}(1-q^{3/2}) (1-q^{7/2})}{(1-q^{5/2})(1+q^5)+\dfrac{q^{5/2}(1-q^{13/2}) (1-q^{17/2})}{(1-q^{5/2}) (1+q^{10})+\ldots}}},
\end{align}
and
\begin{align}\label{S2}
	\nonumber	S_2(q):&=q^{1/4} \dfrac{f(-q^2,-q^8)}{f(-q^3,-q^7)}\\
	&=\dfrac{q^{1/4}(1-q^2)}{(1-q^{5/2})+\dfrac{q^{5/2}(1-q^{1/2}) (1-q^{9/2})}{(1-q^{5/2})(1+q^5)+\dfrac{q^{5/2}(1-q^{11/2}) (1-q^{19/2})}{(1-q^{5/2}) (1+q^{10})+\ldots}}}.
\end{align}
Aricheta and Guadalupe \cite{Aricheta_2024} discussed the modularity of the above two continued fractions and obtained the evaluations as well.
In \cite{Rajkhowa_2026}, Rajkhowa and Saikia studied the following continued fractions of order twenty and derived algebraic relations as well as explicit evaluations.
\begin{align}\label{T1}
	\nonumber	T_1(q):&=q \dfrac{f(-q^3,-q^{17})}{f(-q^7,-q^{13})}\\
	&=\dfrac{q(1-q^3)}{(1-q^5)+\dfrac{q^5(1-q^2) (1-q^8)}{(1-q^5)(1+q^{10})+\dfrac{q^5(1-q^{12}) (1-q^{18})}{(1-q^5) (1+q^{20})+\ldots}}},
\end{align}
and 
\begin{align}\label{T2}
	\nonumber	T_2(q):&=q^{2} \dfrac{f(-q,-q^{19})}{f(-q^9,-q^{11})}\\
	&=\dfrac{q^{2}(1-q)}{(1-q^5)+\dfrac{q^5(1-q^4) (1-q^6)}{(1-q^5)(1+q^{10})+\dfrac{q^5(1-q^{14}) (1-q^{16})}{(1-q^{5}) (1+q^{20})+\ldots}}}.
\end{align}
In \cite{Adiga_2014}, Adiga \textit{et al.} proved identities connecting Ramanujan's cubic continued fraction and a continued fraction of of order six. Further, they obtained several identities associated to these continued fractions. Also, for similar work on continued fractions of order twelve and sixteen, \cite{Adiga_2014b} and \cite{Vanitha_2025} can be referred respectively.

Although numerous identities involving continued fractions $R(q)$, $S_1(q)$, $S_2(q)$, $T_1(q)$ and $T_2(q)$ are known, they are dispersed across the literature. This motivates us to develop a unified approach for deriving identities involving all these continued fractions.  By applying the product expansion of Jacobi's theta function,  the identities relating continued fractions of several different orders can be derived. The purpose of this paper is to obtain a  family of new identities involving Ramanujan's continued fractions of order five, ten and twenty  through the systematic use of Jacobi's theta function $\theta_1(z|\tau)$. By evaluating $\theta_1$ at special arguments $z=\dfrac{k\pi}{20}$ and employing its product representations, we derive explicit expressions for the infinite product 
\begin{align*}
	\Omega_k(q) = \prod_{n=1}^{\infty} (1+\alpha_kq^n +q^{2n} ), \qquad 1\leq k\leq 9,
\end{align*} where the coefficients $\alpha_k$ are algebraic numbers arising from tenth root of unity. These products play an essential role in connecting theta functions to Ramanujan's continued fractions.

This work makes a substantial contribution to the theory of Ramanujan's continued fractions by establishing several new identities involving continued fractions of order twenty. Through the use of Jacobi's theta-function product representations, it uncovers previously unknown relationships between order-twenty continued fractions, order-ten continued fractions, and the Rogers-Ramanujan continued fraction, thereby revealing a richer underlying structure among these classical objects. The work further bridges continued fractions with mock theta functions, Lambert series, and level-twenty theta functions through novel identities derived from logarithmic derivatives and Ramanujan's $_1\psi_1$ summation formula. These results not only provide new analytical tools for studying $q$-series and modular-type functions but also extend and unify earlier developments for continued fractions of orders 6, 12, and 16. 

This research article is classified as follows. In Section 2, we provide all the definitions and preliminary results required. In Section 3, we give the main results that connect all the three continued fractions. In Section 4, we establish Lambert series identities associated to $T_1(q)$ and $T_2(q)$. In Section 5, we derive few more Lambert series identities associated to residue classes modulo twenty.

\section{Essentials}
Let $q \in \mathbb{C}$, such that $ |q| < 1$, throughout. Then $q$-Pochhammer symbol or $q$-shifted factorial \cite{Berndt1991} is defined as follows: For any $\delta \in \mathbb{C}$ and $m \in \mathbb{Z^+}$,
\begin{align*}
	(\delta;q)_0:=1, \quad		(\delta;q)_m:= \prod_{k=1}^{m}(1-\delta q^{k-1}) \quad \text{and} \quad 	(\delta;q)_\infty:= \prod_{k=0}^{\infty}(1-\delta q^{k}).
\end{align*}
\noindent For each $\delta \in \mathbb{C}$, the above product is uniformly and absolutely convergent in every compact subset of the unit disc $|q| < 1$. We use the following notation for ease:
\begin{align*}
	(\delta_1, \delta_2, \ldots, \delta_m;q)_\infty = (\delta_1;q)_\infty (\delta_2;q)_\infty \ldots (\delta_m;q)_\infty.
\end{align*}
The following are the theta functions \cite[ p. 36, Entry 22 (i)-(iii)]{Berndt1991} arising from $f(a,b)$ which are useful in the upcoming sections:
\begin{align}
	\psi(q):=	f(q,q^3)&= \sum_{n=0}^{\infty} q^{\frac{n(n+1)}{2}}=\frac{(q^2;q^2)_\infty}{(q;q^2)_\infty}
\end{align}and 
\begin{align}
	f(-q):= f(-q,-q^2) = \sum_{n=-\infty}^{\infty} (-1)^n  q^{\frac{n(3n-1)}{2}}=(q;q)_\infty.
\end{align}
For $z\ne 0$, the Jacobi's triple product identity states that  
\begin{align*}
	\sum_{n=-\infty}^{\infty} (-1)^n q^{\frac{n(n-1)}{2}} z^n = (q;q)_\infty (z;q)_\infty (q/z;q)_\infty.
\end{align*}
After Ramanujan, 
\begin{align*}
	f(x,y) := (-x;xy)_\infty (-y;xy)_\infty (xy;xy)_\infty, \quad \text{where} \quad |xy|<1.
\end{align*}
Dedekind-eta function is defined by
\begin{align}\label{de}
	\eta(\tau): = q^{1/24} f(-q),
\end{align} where $q=e^{2\pi i \tau}$, $Im \tau >0$.
A Lambert series is a series of the form
\begin{align*}
	L(q) = \sum_{n=1}^{\infty} a_n \dfrac{q^n}{1-q^n}, \qquad |q|<1,
\end{align*} 
where $\{a_n\}_{n\geq 1}$ is an arithmetic sequence. Expanding
\begin{align*}
	\frac{1}{1-q^n}=\sum_{m=0}^{\infty} q^{mn},
\end{align*}
one obtains
\[
\sum_{n=1}^{\infty} \frac{a_n q^n}{1-q^n}
=
\sum_{N=1}^{\infty}
\left(
\sum_{d\mid N} a_d
\right) q^N,
\]
showing that Lambert series naturally generate divisor sums.
The following results are helpful in proving the main results. 

\begin{lemma}[Entry 30 (i), (ii) \& (iii), p.\ 46 of \cite{Adiga_1985}]
	We have 
	\begin{align}
		f(a,ab^2) f(b,a^2b) &= f(a,b) \psi(ab), \label{f1}
	\end{align}
	\begin{align} 
		f(a,b) +f(-a,-b) &= 2f(a^3b,ab^3), \label{f2}
	\end{align}	and 
	\begin{align}
		f(a,b) -f(-a,-b) &= 2af\left(\frac{b}{a},a^5b^3\right). \label{f3}
	\end{align}
\end{lemma}
\noindent Subtracting \eqref{f3} from \eqref{f2}, one can deduce
\begin{align}
	f(-a,-b) =f(a^3b,ab^3) -af\left(\frac{b}{a}, a^5b^3\right). \label{f4}
\end{align}

\section{Identities connecting Ramanujan's continued fractions}
Jacobi's theta function, $\theta_1$, is defined as 
\begin{align}\label{t1}
	\nonumber	\theta_1(z|\tau) &= 2 \sum_{n=0}^{\infty} (-1)^n q^{\frac{(2n+1)^2}{8}} \sin (2n+1) z\\
	&= 2 q^{1/8} \sum_{n=0}^{\infty} (-1)^n q^{\frac{n(n+1)}{2}} \sin (2n+1) z.
\end{align}
\noindent 	In \cite{Chan_2009}, Chan \textit{et al.} showed that 
\begin{align}\label{t2}
	\nonumber	2  \sum_{n=0}^{\infty} (-1)^n q^{\frac{n(n+1)}{2}} \sin (2n+1) z &= \sum_{n=- \infty}^{\infty} (-1)^n q^{\frac{n(n+1)}{2}} \sin (2n+1) z\\
	&= 2 \sin z (q;q)_\infty (qe^{2i z};q)_\infty (qe^{-2i z};q)_\infty
\end{align}
Combining the above two equations, one can see that 
\begin{align}\label{t3}
	\nonumber	\theta_1(z|\tau) &= 2 q^{1/8} \sin z (q;q)_\infty (qe^{2i z};q)_\infty (qe^{-2i z};q)_\infty \\
	&= i q^{1/8} e^{-iz} (q;q)_\infty (e^{2iz};q)_\infty (qe^{-2iz};q)_\infty.
\end{align}
Substituting $z=\frac{\pi}{20}$, $\frac{2\pi}{20}$, $\frac{3\pi}{20}$, $\frac{4\pi}{20}$, $\frac{5\pi}{20}$, $\frac{6\pi}{20}$, $\frac{7\pi}{20}$, $\frac{8\pi}{20}$, and $\frac{9\pi}{20}$ respectively in \eqref{t3}, one can deduce the following:
\begin{align}\label{tk}
	\theta_1\left(\frac{k \pi}{20} | \tau\right) = 2q^{1/8} \left(\sin\frac{k\pi}{20}\right) (q;q)_\infty (qe^{\frac{2k\pi i}{20}};q)_\infty (qe^{\frac{-2k\pi i}{20}};q)_\infty ,
\end{align} for all $k, 1\leq k \leq 9.$
Multiplying all the nine equations together, and using the identities
\begin{align}\label{prodsine}
	\prod_{k=1}^{9} \sin\dfrac{k \pi}{20} = \dfrac{\sqrt{10}}{512}.
\end{align}
and
\begin{align*}
	(1-x) (1-xe^{\frac{20\pi i}{20}})\prod_{k=1}^{9}(1-xe^{\frac{2k\pi i}{20}}) (1-xe^{\frac{-2k\pi i}{20}}) = (1-x^{20}),
\end{align*} one can deduce the following:
\begin{align}\label{tm}
	\prod_{k=1}^{9} \theta_1\left(\frac{k\pi }{20}|\tau\right) = \dfrac{\sqrt{10} \quad \eta^9(\tau) \eta(20\tau)}{\eta(2\tau)}.
\end{align}
Let
\begin{align*}
	&\alpha_1 := -2 \cos \left(\dfrac{2\pi}{20}\right) = -\dfrac{\sqrt{10+2\sqrt{5}}}{2}, \qquad \alpha_6 := -2 \cos \left(\dfrac{12\pi}{20}\right) = \dfrac{\sqrt{5}-1}{2},\\
	&\alpha_2 := -2 \cos \left(\dfrac{4\pi}{20}\right) = -\dfrac{\sqrt{5}+1}{2}, \qquad \qquad \alpha_7 := -2 \cos \left(\dfrac{14\pi}{20}\right) = \dfrac{\sqrt{10-2\sqrt{5}}}{2},\\ 
	&\alpha_3 := -2 \cos \left(\dfrac{6\pi}{20}\right) = -\dfrac{\sqrt{10-2\sqrt{5}}}{2}, \qquad \alpha_8 := -2 \cos \left(\dfrac{16\pi}{20}\right) = \dfrac{\sqrt{5}+1}{2},\\
	&\alpha_4 := -2 \cos \left(\dfrac{8\pi}{20}\right) = -\dfrac{\sqrt{5}-1}{2}, \qquad \qquad \alpha_9 := -2 \cos \left(\dfrac{18\pi}{20}\right) = \dfrac{\sqrt{10+2\sqrt{5}}}{2},\\
	&\alpha_5 := -2 \cos \left(\dfrac{10\pi}{20}\right) =0.
\end{align*}
Using the above and the definition of $\eta(\tau)$, \eqref{tk} can be rewritten as 
\begin{align}\label{Ki}
	\Omega_k(q) := \prod_{n=1}^{\infty} (1+\alpha_kq^n +q^{2n} ) = q^{-\frac{1}{12}} \dfrac{\theta_1(\frac{k \pi}{20}|\tau)}{2\eta(\tau)(\sin\frac{k\pi}{20})},
\end{align} for all $k,1\leq k\leq 9.$
Multiplying the above nine identities, using the identities \eqref{prodsine} and \eqref{tm}, the following identity is obtained.
\begin{align}\label{prodK}
	\prod_{k=1}^{9} \Omega_k(q) = q^{-\frac{3}{4}} \dfrac{\eta(20\tau)}{\eta(2\tau)}.
\end{align} 
We now establish the main results here.
\begin{theorem}
	Let $\Omega_k$, where $1 \leq k \leq 9$, be defined as in \eqref{Ki} where $\alpha_1 = -\dfrac{\sqrt{10+2\sqrt{5}}}{2}$, $\alpha_3 = -\dfrac{\sqrt{10-2\sqrt{5}}}{2}$, $\alpha_7 = \dfrac{\sqrt{10-2\sqrt{5}}}{2}$, and $\alpha_9  = \dfrac{\sqrt{10+2\sqrt{5}}}{2}$. Then, we have the following: 
	\begin{enumerate}[(i)]
		\item \begin{align}\label{O1-O9}
			\nonumber&\Omega_1(q^{1/5}) \prod_{\substack{k=1 \\ k\ne 5}}^{9} \Omega_k(q^{1/5}) - \Omega_9(q^{1/5})  \prod_{\substack{k=1 \\ k\ne 5}}^{9} \Omega_k(q^{1/5}) = - \dfrac{\sqrt{10+2\sqrt{5}} q^{19/40} \eta(4\tau) \psi(-q^5)}{\eta(\tau/5) \eta(4\tau/5)}\\
			\nonumber&+\dfrac{q^{-3/20} \eta(4\tau)\eta(20\tau ) \eta^{1/8}(\tau)}{\eta(\tau/5) \eta(4\tau/5) \eta^{1/8}(5\tau)} \left[  \dfrac{\sqrt{10+2\sqrt{5}}}{\sqrt{T_1(q)} \sqrt[4]{S_2(q)} \sqrt[8]{R(q)}}
			\right. \\ & \left. +\dfrac{\sqrt{50-10\sqrt{5}}+3\sqrt{10-2\sqrt{5}}}{2} \left\{\dfrac{\sqrt{T_1(q)}}{\sqrt[4]{S_2(q)} \sqrt[8]{R(q)}}-\sqrt{T_1(q)} \sqrt[4]{S_1(q)} \sqrt[8]{R(q)} \right\}   \right],
		\end{align}\\
		\item 	\begin{align}\label{O3-O7}
			\nonumber&\Omega_3(q^{1/5}) \prod_{\substack{k=1 \\ k\ne 5}}^{9} \Omega_k(q^{1/5}) - \Omega_7(q^{1/5})  \prod_{\substack{k=1 \\ k\ne 5}}^{9} \Omega_k(q^{1/5}) 
			= \dfrac{\sqrt{10+2\sqrt{5}} (\sqrt{5}-1)q^{19/40} \eta(4\tau) \psi(-q^5)}{2\eta(\tau/5) \eta(4\tau/5)}\\
			\nonumber&-\dfrac{q^{-3/20} \eta(4\tau)\eta(20\tau ) \eta^{1/8}(\tau)}{\eta(\tau/5) \eta(4\tau/5) \eta^{1/8}(5\tau)} \left[  \dfrac{\sqrt{10+2\sqrt{5}} (\sqrt{5}-1)}{2\sqrt{T_1(q)} \sqrt[4]{S_2(q)} \sqrt[8]{R(q)}} \right. \\ & \left. + \dfrac{\sqrt{10-2\sqrt{5}} (\sqrt{5}-1)}{2} \left\{\sqrt{T_1(q)} \sqrt[4]{S_1(q)} \sqrt[8]{R(q)}-\dfrac{\sqrt{T_1(q)}}{\sqrt[4]{S_2(q)} \sqrt[8]{R(q)}} \right\} \right],
		\end{align}\\
		\item 	\begin{align}\label{1O1+9O9}
			\nonumber&(1+\alpha_9)\Omega_9(q^{1/5}) \prod_{\substack{k=1 \\ k\ne 5}}^{9} \Omega_k(q^{1/5}) + (1+\alpha_1)\Omega_1(q^{1/5})  \prod_{\substack{k=1 \\ k\ne 5}}^{9} \Omega_k(q^{1/5}) =-\dfrac{2q^{19/40} \eta(4\tau) \psi(-q^5)}{\eta(\tau/5) \eta(4\tau/5)}\\
			\nonumber&+\dfrac{q^{-3/20} \eta(4\tau)\eta(20\tau ) \eta^{1/8}(\tau)}{\eta(\tau/5) \eta(4\tau/5) \eta^{1/8}(5\tau)}
			\left[  2 \dfrac{\sqrt[4]{S_1(q)} \sqrt[8]{R(q)}}{\sqrt{T_2(q)}}+ 2(\sqrt{5}+1) \dfrac{\sqrt{T_1(q)}}{\sqrt[4]{S_2(q)} \sqrt[8]{R(q)}}
			\right. \\ & \left.
			+(3+\sqrt{5}) \left\{\dfrac{1}{\sqrt{T_1(q)} \sqrt[4]{S_2(q)} \sqrt[8]{R(q)}}-\sqrt{T_1(q)} \sqrt[4]{S_1(q)} \sqrt[8]{R(q)}\right\} \right],
		\end{align}\\
		\item 	\begin{align}\label{1O1-9O9}
			(1+\alpha_9)\Omega_9(q^{1/5}) \prod_{\substack{k=1 \\ k\ne 5}}^{9} \Omega_k(q^{1/5}) - (1+\alpha_1)\Omega_1(q^{1/5})  \prod_{\substack{k=1 \\ k\ne 5}}^{9} \Omega_k(q^{1/5}) = 0,
		\end{align}\\
		\item 	\begin{align}\label{7O7+3O3}
			\nonumber&(1+\alpha_7)\Omega_7(q^{1/5}) \prod_{\substack{k=1 \\ k\ne 5}}^{9} \Omega_k(q^{1/5}) + (1+\alpha_3)\Omega_3(q^{1/5})  \prod_{\substack{k=1 \\ k\ne 5}}^{9} \Omega_k(q^{1/5}) 			
			= -\dfrac{2q^{19/40} \eta(4\tau) \psi(-q^5)}{\eta(\tau/5) \eta(4\tau/5)}\\
			\nonumber&
			+\dfrac{q^{-3/20} \eta(4\tau)\eta(20\tau ) \eta^{1/8}(\tau)}{\eta(\tau/5) \eta(4\tau/5) \eta^{1/8}(5\tau)} \left[ \dfrac{2\sqrt[4]{S_1(q)} \sqrt[8]{R(q)}}{\sqrt{T_2(q)}} -\dfrac{2(\sqrt{5}-1) \sqrt{T_1(q)}}{\sqrt[4]{S_2(q)} \sqrt[8]{R(q)}}
			\right. \\ & \left.
			(\sqrt{5}-3) \left\{\sqrt{T_1(q)} \sqrt[4]{S_1(q)} \sqrt[8]{R(q)}-\dfrac{1}{\sqrt{T_1(q)} \sqrt[4]{S_2(q)} \sqrt[8]{R(q)}}\right\}\right] ,
		\end{align}\\
		\item 	\begin{align}\label{7O7-3O3}
			\nonumber&(1+\alpha_7)\Omega_7(q^{1/5}) \prod_{\substack{k=1 \\ k\ne 5}}^{9} \Omega_k(q^{1/5}) - (1+\alpha_3)\Omega_3(q^{1/5})  \prod_{\substack{k=1 \\ k\ne 5}}^{9} \Omega_k(q^{1/5})\\
			\nonumber&
			= \dfrac{ \left(\sqrt{10-2\sqrt{5}}-\sqrt{10+2\sqrt{5}}\right)q^{19/40} \eta(4\tau) \psi(-q^5)}{\eta(\tau/5) \eta(4\tau/5)}
			+\dfrac{q^{-3/20} \eta(4\tau)\eta(20\tau ) \eta^{1/8}(\tau)}{\eta(\tau/5) \eta(4\tau/5) \eta^{1/8}(5\tau)} \\
			\nonumber&
			\left[ \sqrt{10-2\sqrt{5}} \left\{ \dfrac{\sqrt[4]{S_1(q)} \sqrt[8]{R(q)}}{\sqrt{T_2(q)}} - \dfrac{\sqrt{T_1(q)}}{\sqrt[4]{S_2(q)} \sqrt[8]{R(q)}} \right\} \right. \\ & \left.+ \left(\sqrt{10-2\sqrt{5}}-\sqrt{10+2\sqrt{5}}\right) \sqrt{T_1(q)} \sqrt[4]{S_1(q)} \sqrt[8]{R(q)}     \right],
		\end{align}\\
		\item 	\begin{align}\label{O99+O11}
			\nonumber& \alpha_9 \Omega_9(q^{1/5}) \prod_{\substack{k=1 \\ k\ne 5}}^{9} \Omega_k(q^{1/5}) + \alpha_1  \Omega_1(q^{1/5})  \prod_{\substack{k=1 \\ k\ne 5}}^{9} \Omega_k(q^{1/5}) = -\dfrac{ (5+\sqrt{5}) q^{19/40} \eta(4\tau) \psi(-q^5)}{\eta(\tau/5) \eta(4\tau/5)}\\
			\nonumber&+\dfrac{q^{-3/20} \eta(4\tau)\eta(20\tau ) \eta^{1/8}(\tau)}{\eta(\tau/5) \eta(4\tau/5) \eta^{1/8}(5\tau)} \left[   (5+3\sqrt{5})\left\{\dfrac{\sqrt{T_1(q)}}{ \sqrt[4]{S_2(q)} \sqrt[8]{R(q)}}-\sqrt{T_1(q)}\sqrt[4]{S_1(q)} \sqrt[8]{R(q)} \right\}
			\right. \\ & \left. +
			\dfrac{ (5+\sqrt{5})}{\sqrt{T_1(q)} \sqrt[4]{S_2(q)} \sqrt[8]{R(q)}}
			\right],
		\end{align}\\
		\item	\begin{align}\label{O99-O11}
			\nonumber& \alpha_9 \Omega_9(q^{1/5}) \prod_{\substack{k=1 \\ k\ne 5}}^{9} \Omega_k(q^{1/5}) - \alpha_1  \Omega_1(q^{1/5})  \prod_{\substack{k=1 \\ k\ne 5}}^{9} \Omega_k(q^{1/5})\\
			\nonumber &
			=\dfrac{\left(\sqrt{10-2\sqrt{5}}+2\sqrt{10+2\sqrt{5}}\right) q^{19/40}  \eta(4\tau) \psi(-q^5)}{\eta(\tau/5) \eta(4\tau/5)}
			+\dfrac{q^{-3/20} \eta(4\tau)\eta(20\tau ) \eta^{1/8}(\tau)}{\eta(\tau/5) \eta(4\tau/5) \eta^{1/8}(5\tau)} \\
			\nonumber&  \left[ \sqrt{10+2\sqrt{5}} \left\{ \dfrac{\sqrt[4]{S_1(q)} \sqrt[8]{R(q)}}{\sqrt{T_2(q)}} - \dfrac{1}{\sqrt{T_1(q)}\sqrt[4]{S_2(q)} \sqrt[8]{R(q)}} \right\} \right. \\ \nonumber& \left.- \left(\sqrt{10-2\sqrt{5}}+2\sqrt{10+2\sqrt{5}}\right) \dfrac{\sqrt{T_1(q)}}{ \sqrt[4]{S_2(q)} \sqrt[8]{R(q)} }  
			\right. \\ & \left.+ \left(2\sqrt{10-2\sqrt{5}}+2\sqrt{10+2\sqrt{5}}\right) \sqrt{T_1(q)} \sqrt[4]{S_2(q)} \sqrt[8]{R(q)}  \right],
		\end{align}\\
		\item	\begin{align}\label{O77+O33}
			\nonumber& \alpha_7 \Omega_7(q^{1/5}) \prod_{\substack{k=1 \\ k\ne 5}}^{9} \Omega_k(q^{1/5}) + \alpha_3  \Omega_3(q^{1/5})  \prod_{\substack{k=1 \\ k\ne 5}}^{9} \Omega_k(q^{1/5}) = -\dfrac{ (5-\sqrt{5}) q^{19/40} \eta(4\tau) \psi(-q^5)}{\eta(\tau/5) \eta(4\tau/5)}\\
			\nonumber&+\dfrac{q^{-3/20} \eta(4\tau)\eta(20\tau ) \eta^{1/8}(\tau)}{\eta(\tau/5) \eta(4\tau/5) \eta^{1/8}(5\tau)} \left[   (3\sqrt{5}-5)\left\{\sqrt{T_1(q)}\sqrt[4]{S_1(q)} \sqrt[8]{R(q)}-\dfrac{\sqrt{T_1(q)}}{ \sqrt[4]{S_2(q)} \sqrt[8]{R(q)}} \right\}
			\right. \\ & \left. +
			\dfrac{ (5-\sqrt{5})}{\sqrt{T_1(q)} \sqrt[4]{S_2(q)} \sqrt[8]{R(q)}}
			\right],
		\end{align} and\\
		\item	\begin{align}\label{O77-O33}
			\nonumber& \alpha_7 \Omega_7(q^{1/5}) \prod_{\substack{k=1 \\ k\ne 5}}^{9} \Omega_k(q^{1/5}) - \alpha_3  \Omega_3(q^{1/5})  \prod_{\substack{k=1 \\ k\ne 5}}^{9} \Omega_k(q^{1/5})\\
			\nonumber&= \dfrac{q^{19/40} \left(2\sqrt{10-2\sqrt{5}}-\sqrt{10+2\sqrt{5}}\right) \eta(4\tau) \psi(-q^5)}{\eta(\tau/5) \eta(4\tau/5)}
			+\dfrac{q^{-3/20} \eta(4\tau)\eta(20\tau ) \eta^{1/8}(\tau)}{\eta(\tau/5) \eta(4\tau/5) \eta^{1/8}(5\tau)} \\
			\nonumber&  \left[ \sqrt{10-2\sqrt{5}} \left\{ \dfrac{\sqrt[4]{S_1(q)} \sqrt[8]{R(q)}}{\sqrt{T_2(q)}} - \dfrac{1}{\sqrt{T_1(q)}\sqrt[4]{S_2(q)} \sqrt[8]{R(q)}} \right\} \right. \\ \nonumber& \left.- \left(2\sqrt{10-2\sqrt{5}}-\sqrt{10+2\sqrt{5}}\right) \dfrac{\sqrt{T_1(q)}}{ \sqrt[4]{S_2(q)} \sqrt[8]{R(q)} }  
			\right. \\ & \left.+ \left(2\sqrt{10-2\sqrt{5}}-2\sqrt{10+2\sqrt{5}}\right) \sqrt{T_1(q)} \sqrt[4]{S_2(q)} \sqrt[8]{R(q)}  \right],
		\end{align}
	\end{enumerate}
	where $R(q)$ is the Rogers-Ramanujan continued fraction, $S_1(q)$ and $S_2(q)$ are the continued fractions of order ten given by \eqref{S1} and \eqref{S2} and $T_1(q)$ and $T_2(q)$ are the continued fractions of order twenty given by \eqref{T1} and \eqref{T2} respectively.
\end{theorem}
\begin{proof}[Proof of \eqref{O1-O9}]
	Consider 
	\begin{align*}
		\Omega_1(q) - \Omega_9(q) &= \prod_{n=1}^{\infty} (1+\alpha_1 q^n +q^{2n}) - \prod_{n=1}^{\infty} (1+\alpha_9 q^n +q^{2n})\\
		&= \dfrac{q^{-1/12}}{\eta(\tau)} \left(\dfrac{\theta_1\left(\frac{\pi}{20}|\tau\right)}{2 \sin \left(\frac{\pi}{20}\right)}-\dfrac{\theta_1\left(\frac{9\pi}{20}|\tau\right)}{2 \sin \left(\frac{9\pi}{20}\right)}\right).
	\end{align*}
	Using the identity \eqref{t1}, the above identity can be written as
	\begin{align}\label{3.1.0}
		\Omega_1(q) - \Omega_9(q) 
		&= \dfrac{q^{-1/12}}{\eta(\tau)} \sum_{n=0}^{\infty} (-1)^n A_1(n) q^{\frac{(2n+1)^2}{8}},
	\end{align} where
	\begin{align*}
		A_1(n)= \dfrac{\sin (2n+1) \frac{\pi}{20}}{\sin \frac{\pi}{20}}- \dfrac{\sin (2n+1) \frac{9\pi}{20}}{\sin \frac{9\pi}{20}}.
	\end{align*}
	Using MAPLE computations, one can see that
	\begin{align*}
		&A_1(20m+0) =0,  \\
		& A_1(20m+1)=\sqrt{10+2\sqrt{5}}, \qquad A_1(20m+2)=\sqrt{10+2\sqrt{5}}, \\
		&A_1(20m+3)=\dfrac{\sqrt{50-10\sqrt{5}}+3\sqrt{10-2\sqrt{5}}}{2}, \\
		&A_1(20m+4)=\dfrac{\sqrt{50-10\sqrt{5}}+3\sqrt{10-2\sqrt{5}}}{2}, \\
		&A_1(20m+5)=\dfrac{\sqrt{50-10\sqrt{5}}+3\sqrt{10-2\sqrt{5}}}{2}, \\
		&A_1(20m+6)=\dfrac{\sqrt{50-10\sqrt{5}}+3\sqrt{10-2\sqrt{5}}}{2}, \\
		&A_1(20m+7)=\sqrt{10+2\sqrt{5}}, \qquad
		A_1(20m+8)=\sqrt{10+2\sqrt{5}}, \\
		&A_1(20m+9) =0, \quad \qquad \qquad \qquad A_1(20m+10) = 0,\\
		&A_1(20m+11)=-\sqrt{10+2\sqrt{5}}, 	\quad	A_1(20m+12)=-\sqrt{10+2\sqrt{5}},\\
		&A_1(20m+13)=-\dfrac{\sqrt{50-10\sqrt{5}}+3\sqrt{10-2\sqrt{5}}}{2},\\
		&A_1(20m+14)=-\dfrac{\sqrt{50-10\sqrt{5}}+3\sqrt{10-2\sqrt{5}}}{2},\\
		&A_1(20m+15)=-\dfrac{\sqrt{50-10\sqrt{5}}+3\sqrt{10-2\sqrt{5}}}{2},\\
		&A_1(20m+16)=-\dfrac{\sqrt{50-10\sqrt{5}}+3\sqrt{10-2\sqrt{5}}}{2},\\
		&A_1(20m+17)=-\sqrt{10+2\sqrt{5}},\quad		A_1(20m+18)=-\sqrt{10+2\sqrt{5}},\\
		&A_1(20m+19) = 0.
	\end{align*}
	Thus,
	\begin{align*}
		&\sum_{n=0}^{\infty} (-1)^n A_1(n) q^{\frac{(2n+1)^2}{8}} =  \left(\sqrt{10+2\sqrt{5}}\right)\left\{-\sum_{m=0}^{\infty} q^{\frac{(40m+3)^2}{8}} 
		+\sum_{m=0}^{\infty}  q^{\frac{(40m+5)^2}{8}}\right\}\\
		&+\left(\dfrac{\sqrt{50-10\sqrt{5}}+3\sqrt{10-2\sqrt{5}}}{2}\right)\left\{-\sum_{m=0}^{\infty} q^{\frac{(40m+7)^2}{8}}
		+\sum_{m=0}^{\infty} q^{\frac{(40m+9)^2}{8}} \right. \\
		& \left.-\sum_{m=0}^{\infty}  q^{\frac{(40m+11)^2}{8}}
		+\sum_{m=0}^{\infty}  q^{\frac{(40m+13)^2}{8}}\right\} + \left(\sqrt{10+2\sqrt{5}}\right) \left\{-\sum_{m=0}^{\infty} q^{\frac{(40m+15)^2}{8}}\right. \\
		& \left.
		+\sum_{m=0}^{\infty}  q^{\frac{(40m+17)^2}{8}}
		+\sum_{m=0}^{\infty}   q^{\frac{(40m+23)^2}{8}} -\sum_{m=0}^{\infty}   q^{\frac{(40m+25)^2}{8}}\right\} 
		+\left(\dfrac{\sqrt{50-10\sqrt{5}}+3\sqrt{10-2\sqrt{5}}}{2}\right)\\ & \times \left\{ \sum_{m=0}^{\infty}  q^{\frac{(40m+27)^2}{8}}
		+\sum_{m=0}^{\infty}   q^{\frac{(40m+29)^2}{8}}
		+\sum_{m=0}^{\infty}   q^{\frac{(40m+31)^2}{8}}
		+\sum_{m=0}^{\infty}   q^{\frac{(40m+33)^2}{8}}\right\}
		+\left(\sqrt{10+2\sqrt{5}}\right)\\ & \times \left\{\sum_{m=0}^{\infty} q^{\frac{(40m+35)^2}{8}}
		-\sum_{m=0}^{\infty}  q^{\frac{(40m+37)^2}{8}}\right\}.
	\end{align*}
	Changing $m$ to $-m$ and $m$ to $m-1$ in the first and last eight summations respectively, above equation deduces to 
	
	\begin{align*}
		&\sum_{n=0}^{\infty} (-1)^n A_1(n) q^{\frac{(2n+1)^2}{8}} =  \left(\sqrt{10+2\sqrt{5}}\right)\left\{-\sum_{m=-\infty}^{\infty} q^{\frac{(40m-3)^2}{8}} 
		+\sum_{m=-\infty}^{\infty}  q^{\frac{(40m-5)^2}{8}}\right. \\
		& \left.-\sum_{m=-\infty}^{\infty} q^{\frac{(40m-15)^2}{8}}
		+\sum_{m=-\infty}^{\infty}  q^{\frac{(40m-17)^2}{8}}\right\} +\left(\dfrac{\sqrt{50-10\sqrt{5}}+3\sqrt{10-2\sqrt{5}}}{2}\right)\\ &\times \left\{-\sum_{m=-\infty}^{\infty} q^{\frac{(40m-7)^2}{8}}
		+\sum_{m=-\infty}^{\infty} q^{\frac{(40m-9)^2}{8}} -\sum_{m=-\infty}^{\infty}  q^{\frac{(40m-11)^2}{8}}
		+\sum_{m=-\infty}^{\infty}  q^{\frac{(40m-13)^2}{8}}\right\} .
	\end{align*}
	Using the definition of $f(a,b)$ in the above, the above expression reduces to 
	\begin{align}\label{3.1.1}
		\nonumber	&\sum_{n=0}^{\infty} (-1)^n A_1(n) q^{\frac{(2n+1)^2}{8}} = -\left(\sqrt{10+2\sqrt{5}}\right) q^{9/8} \left\{ f(q^{170},q^{230})-q^2 f(q^{150}, q^{250}) 
		\right. \\ & \left. \nonumber -q^{35} f(q^{30}, q^{370}) +q^{27} f(q^{50}, q^{350}) \right\} +\left(\dfrac{\sqrt{50-10\sqrt{5}}+3\sqrt{10-2\sqrt{5}}}{2}\right) q^{9/8} 
		\\ &  \left\{ q^5 f(q^{130}, q^{270})
		-q^9 f(q^{110},q^{290}) +q^{14} f(q^{90}, q^{310})-q^{20} f(q^{70}, q^{330})\right\}.
	\end{align}
	Setting $(a,b)=(q^5,q^{95}),(q^{15}, q^{85}), (q^{25},q^{75})$ and $(q^{35},q^{65})$ in \eqref{f4} yields respectively the following equations:
	\begin{align}\label{fab1}
		f(-q^5,-q^{95}) &= f(q^{110}, q^{290}) -q^5 f(q^{90}, q^{310}),
	\end{align}
	\begin{align}\label{fab2}
		f(-q^{15},-q^{85}) &= f(q^{130}, q^{270}) -q^{15} f(q^{70}, q^{330}),
	\end{align}
	\begin{align}\label{fab3}
		f(-q^{25},-q^{75}) &= f(q^{150}, q^{250}) -q^{25} f(q^{50}, q^{350}),
	\end{align}and
	\begin{align}\label{fab4}
		f(-q^{35},-q^{65}) &= f(q^{170}, q^{230}) -q^{35} f(q^{30}, q^{370}).
	\end{align}
	Using these in \eqref{3.1.1}, 
	\begin{align}\label{3.1.2}
		\nonumber	&\sum_{n=0}^{\infty} (-1)^n A_1(n) q^{\frac{(2n+1)^2}{8}} = -\left(\sqrt{10+2\sqrt{5}}\right) q^{9/8} \left\{f(-q^{35},-q^{65}) -q^2f(-q^{25}, -q^{})\right\}\\
		&+\left(\dfrac{\sqrt{50-10\sqrt{5}}+3\sqrt{10-2\sqrt{5}}}{2}\right) q^{9/8}\left\{q^5 f(-q^{15},-q^{85}) -q^9 f(-q^5,-q^{95})\right\}.
	\end{align}
	Utilizing \eqref{3.1.2} in \eqref{3.1.0}, one can see that
	\begin{align}\label{3.1.3}
		\nonumber &	\Omega_1(q) - \Omega_9(q) 
		= \dfrac{q^{25/24}}{\eta(\tau)} \left[\left(-\sqrt{10+2\sqrt{5}}\right)   \left\{f(-q^{35},-q^{65}) -q^2f(-q^{25}, -q^{75})\right\} \right. \\
		& \left.+\left(\dfrac{\sqrt{50-10\sqrt{5}}+3\sqrt{10-2\sqrt{5}}}{2}\right)  \left\{q^5 f(-q^{15},-q^{85}) -q^9 f(-q^5,-q^{95})\right\}\right].
	\end{align}
	Multiplying on both the sides of \eqref{3.1.3} by \eqref{prodK},
	\begin{align}\label{3.1.4}
		\nonumber &	\Omega_1(q)\prod_{k=1}^{9} \Omega_k - \Omega_9(q) \Omega_1(q)\prod_{k=1}^{9} \Omega_k
		= \dfrac{q^{7/24}\eta(20\tau)}{\eta(\tau)\eta(2\tau)}\\
		& \times \nonumber \left[ \left(-\sqrt{10+2\sqrt{5}}\right)   \left\{f(-q^{35},-q^{65}) -q^2f(-q^{25}, -q^{75})\right\} \right.\\
		& \left.+\left(\dfrac{\sqrt{50-10\sqrt{5}}+3\sqrt{10-2\sqrt{5}}}{2}\right)  \left\{q^5 f(-q^{15},-q^{85}) -q^9 f(-q^5,-q^{95})\right\}\right].
	\end{align}
	Noting the fact that $\alpha_5=0$ and $(-q^2;q^2)_\infty = q^{-2/24}\dfrac{\eta(4\tau)}{\eta(2\tau)}$, upon simplification and changing $q$ to $q^{1/5}$ throughout \eqref{3.1.4}, the identity \eqref{O1-O9} is obtained.
	The proof of the identity \eqref{O3-O7} is similar to the above and hence is omitted. 
\end{proof}
\begin{proof}[Proof of \eqref{1O1+9O9}]
	Consider 
	\begin{align*}
		&(1+\alpha_9)	\Omega_9(q) + (1+\alpha_1) \Omega_1(q)\\ &=(1+\alpha_9)\prod_{n=1}^{\infty} (1+\alpha_9 q^n +q^{2n}) + (1+\alpha_1)\prod_{n=1}^{\infty} (1+\alpha_1 q^n +q^{2n})\\
		&= \dfrac{q^{-1/12}}{\eta(\tau)} \left(\dfrac{(1+\alpha_9)\theta_1\left(\frac{9\pi}{20}|\tau\right)}{2 \sin \left(\frac{9\pi}{20}\right)}+\dfrac{(1+\alpha_1)\theta_1\left(\frac{\pi}{20}|\tau\right)}{2 \sin \left(\frac{\pi}{20}\right)}\right).
	\end{align*}
	Using the identity \eqref{t1}, the above identity can be written as
	\begin{align}\label{3.1.5}
		(1+\alpha_9)	\Omega_9(q) + (1+\alpha_1) \Omega_1(q)
		&= \dfrac{q^{-1/12}}{\eta(\tau)} \sum_{n=0}^{\infty} (-1)^n A_2(n) q^{\frac{(2n+1)^2}{8}},
	\end{align} where
	\begin{align*}
		A_2(n)= \dfrac{(1+\alpha_9)\sin (2n+1) \frac{9\pi}{20}}{\sin \frac{9\pi}{20}}+ \dfrac{(1+\alpha_1)\sin (2n+1) \frac{\pi}{20}}{\sin \frac{\pi}{20}}.
	\end{align*}
	Using MAPLE computations, one can see that
	\begin{align*}
		A_2(20m+0) =2,  &\qquad A_2(20m+10) = -2,\\
		A_2(20m+1)=-(3+\sqrt{5}), &\qquad A_2(20m+11)=(3+\sqrt{5}),\\
		A_2(20m+2)=-2, &\qquad A_2(20m+12)=2,\\
		A_2(20m+3)=-2(\sqrt{5}+1), &\qquad A_2(20m+13)=2(\sqrt{5}+1) ,\\
		A_2(20m+4)=-(3+\sqrt{5}), &\qquad A_2(20m+14)=(3+\sqrt{5}),\\
		A_2(20m+5)=-(3+\sqrt{5}), &\qquad A_2(20m+15)=(3+\sqrt{5}), \\
		A_2(20m+6)=-2(\sqrt{5}+1), &\qquad A_2(20m+16)=2(\sqrt{5}+1),\\
		A_2(20m+7)=-2, &\qquad A_2(20m+17)=2, \\
		A_2(20m+8)=-(3+\sqrt{5}), &\qquad A_2(20m+18)=  (3+\sqrt{5}),\\
		A_2(20m+9) =2, &\qquad A_2(20m+19) =-2.  
	\end{align*}
	Thus,
	\begin{align*}
		&\sum_{n=0}^{\infty} (-1)^n A_2(n) q^{\frac{(2n+1)^2}{8}} = 2\sum_{m=0}^{\infty} q^{\frac{(40m+1)^2}{8}}+ (3+\sqrt{5})\sum_{m=0}^{\infty} q^{\frac{(40m+3)^2}{8}} 
		-2\sum_{m=0}^{\infty} q^{\frac{(40m+5)^2}{8}}\\
		&+2(\sqrt{5}+1)\sum_{m=0}^{\infty} q^{\frac{(40m+7)^2}{8}}
		-(3+\sqrt{5})\sum_{m=0}^{\infty} q^{\frac{(40m+9)^2}{8}}
		+(3+\sqrt{5})\sum_{m=0}^{\infty}q^{\frac{(40m+11)^2}{8}}
		\\
		&-2(\sqrt{5}+1)\sum_{m=0}^{\infty} q^{\frac{(40m+13)^2}{8}}+2\sum_{m=0}^{\infty}  q^{\frac{(40m+15)^2}{8}}
		-(3+\sqrt{5})\sum_{m=0}^{\infty}  q^{\frac{(40m+17)^2}{8}}\\
		&
		-2\sum_{m=0}^{\infty}  q^{\frac{(40m+19)^2}{8}}-2\sum_{m=0}^{\infty}  q^{\frac{(40m+21)^2}{8}}-(3+\sqrt{5})\sum_{m=0}^{\infty}  q^{\frac{(40m+23)^2}{8}}+2\sum_{m=0}^{\infty}  q^{\frac{(40m+25)^2}{8}}\\
		&
		-2(\sqrt{5}+1)\sum_{m=0}^{\infty}  q^{\frac{(40m+27)^2}{8}} +(3+\sqrt{5})\sum_{m=0}^{\infty}  q^{\frac{(40m+29)^2}{8}}
		-(3+\sqrt{5})\sum_{m=0}^{\infty}  q^{\frac{(40m+31)^2}{8}}\\
		&+2(\sqrt{5}+1)\sum_{m=0}^{\infty}   q^{\frac{(40m+33)^2}{8}}
		-2\sum_{m=0}^{\infty}   q^{\frac{(40m+35)^2}{8}}+(3+\sqrt{5})\sum_{m=0}^{\infty}   q^{\frac{(40m+37)^2}{8}}+2\sum_{m=0}^{\infty}   q^{\frac{(40m+39)^2}{8}}.
	\end{align*}
	Changing $m$ to $-m$ and $m$ to $m-1$ in the first and last ten summations respectively, above equation deduces to 
	
	\begin{align*}
		&\sum_{n=0}^{\infty} (-1)^n A_2(n) q^{\frac{(2n+1)^2}{8}} = 2(\sqrt{5}+1) \left\{\sum_{m=-\infty}^{\infty} q^{\frac{(40m-7)^2}{8}}
		-\sum_{m=-\infty}^{\infty}   q^{\frac{(40m-13)^2}{8}}\right\}\\ &+2\left\{\sum_{m=-\infty}^{\infty} q^{\frac{(40m-1)^2}{8}} 
		-\sum_{m=-\infty}^{\infty}  q^{\frac{(40m-5)^2}{8}}+\sum_{m=-\infty}^{\infty}  q^{\frac{(40m-15)^2}{8}} 
		-\sum_{m=-\infty}^{\infty}  q^{\frac{(40m-19)^2}{8}}\right\}\\ &			+(3+\sqrt{5})\left\{
		\sum_{m=-\infty}^{\infty}  q^{\frac{(40m-3)^2}{8}}
		-\sum_{m=-\infty}^{\infty}  q^{\frac{(40m-9)^2}{8}}+\sum_{m=-\infty}^{\infty}   q^{\frac{(40m-11)^2}{8}} 
		-\sum_{m=-\infty}^{\infty}   q^{\frac{(40m-17)^2}{8}}\right\} .
	\end{align*}
	Using the definition of $f(a,b)$, the above expression reduces to 
	\begin{align}\label{3.1.6}
		\nonumber	&\sum_{n=0}^{\infty} (-1)^n A_2(n) q^{\frac{(2n+1)^2}{8}} = 2(\sqrt{5}+1) q^{1/8}\left\{ q^{6} f(q^{130}, q^{270}) -q^{21} f(q^{70}, q^{330}) \right\}
		\\ &  \nonumber +2q^{1/8}\left\{f(q^{190},q^{210})-q^3 f(q^{150}, q^{250}) +q^{28}f(q^{50},q^{350}) -q^{45} f(q^{10},q^{390})\right\}\\ & +(3+\sqrt{5}) q^{1/8} 
		\left\{q f(q^{170}, q^{230}) 
		-q^{10} f(q^{110}, q^{290}) +q^{15} f(q^{90},q^{310}) -q^{36} f(q^{30}, q^{370})\right\}.
	\end{align}
	Setting $(a,b)=(q^{45},q^{55})$ in \eqref{f4} yields   the following equation. 
	\begin{align}\label{fab5}
		f(-q^{45},-q^{55}) &= f(q^{190}, q^{210}) -q^{45} f(q^{10}, q^{390}).
	\end{align}
	Using \eqref{fab1}, \eqref{fab2}, \eqref{fab3}, \eqref{fab4} and \eqref{fab5}, \eqref{3.1.6} becomes 
	\begin{align}\label{3.1.7}
		\nonumber	\sum_{n=0}^{\infty} (-1)^n &A_2(n) q^{\frac{(2n+1)^2}{8}} = 2(\sqrt{5}+1) q^{49/8}f(-q^{15},-q^{85})\\ &\nonumber+2q^{1/8}\left\{f(-q^{45},-q^{55}) -q^3f(-q^{25}, -q^{75})\right\}\\
		&+(3+\sqrt{5}) q^{1/8} \left\{qf(-q^{35},-q^{65}) -q^{10} f(-q^5,-q^{95})\right\}.
	\end{align}
	Utilizing \eqref{3.1.7} in \eqref{3.1.5}, one can see that
	\begin{align}\label{3.1.8}
		\nonumber 	(1+\alpha_9)	\Omega_9(q) &+ (1+\alpha_1) \Omega_1(q)
		= \dfrac{q^{-1/12}}{\eta(\tau)} \left[2(\sqrt{5}+1) q^{49/8}f(-q^{15},-q^{85}) \right. \\ & \left. \nonumber+2q^{1/8}\left\{f(-q^{45},-q^{55}) -q^3f(-q^{25}, -q^{75})\right\} \right. \\ & \left.+(3+\sqrt{5}) q^{1/8} \left\{qf(-q^{35},-q^{65}) -q^{10} f(-q^5,-q^{95})\right\}\right].
	\end{align}
	Multiplying on both the sides of \eqref{3.1.8} by \eqref{prodK},
	\begin{align}\label{3.1.9}
		\nonumber &	(1+\alpha_9)\Omega_9(q)\prod_{k=1}^{9} \Omega_k + (1+\alpha_1) \Omega_1(q) \Omega_1(q)\prod_{k=1}^{9} \Omega_k
		= \dfrac{q^{-5/6}\eta(20\tau)}{\eta(\tau)\eta(2\tau)}\\ & \nonumber \times \left[2(\sqrt{5}+1) q^{49/8}f(-q^{15},-q^{85}) +2q^{1/8}\left\{f(-q^{45},-q^{55}) -q^3f(-q^{25}, -q^{75})\right\} \right. \\ & \left.+(3+\sqrt{5}) q^{1/8} \left\{qf(-q^{35},-q^{65}) -q^{10} f(-q^5,-q^{95})\right\}\right].
	\end{align}
	Noting the fact that $\alpha_5=0$ and $(-q^2;q^2)_\infty = q^{-2/24}\dfrac{\eta(4\tau)}{\eta(2\tau)}$, upon simplification and changing $q$ to $q^{1/5}$ throughout \eqref{3.1.9}, the identity \eqref{1O1+9O9} is obtained. The proofs of the identities \eqref{1O1-9O9}, \eqref{7O7+3O3} and \eqref{7O7-3O3} are similar to the above and hence are omitted. 
\end{proof}
\begin{proof}[Proof of \eqref{O99-O11}]
	Consider 
	\begin{align*}
		\alpha_9	\Omega_9(q) - \alpha_1 \Omega_1(q) &=\alpha_9\prod_{n=1}^{\infty} (1+\alpha_9 q^n +q^{2n}) - \alpha_1\prod_{n=1}^{\infty} (1+\alpha_1 q^n +q^{2n})\\
		&= \dfrac{q^{-1/12}}{\eta(\tau)} \left(\dfrac{\alpha_9\theta_1\left(\frac{9\pi}{20}|\tau\right)}{2 \sin \left(\frac{9\pi}{20}\right)}+\dfrac{\alpha_1\theta_1\left(\frac{\pi}{20}|\tau\right)}{2 \sin \left(\frac{\pi}{20}\right)}\right).
	\end{align*}
	Using the identity \eqref{t1}, the above identity can be written as
	\begin{align}\label{3.1.10}
		\alpha_9	\Omega_9(q) - \alpha_1 \Omega_1(q)
		&= \dfrac{q^{-1/12}}{\eta(\tau)} \sum_{n=0}^{\infty} (-1)^n A_3(n) q^{\frac{(2n+1)^2}{8}},
	\end{align} where
	\begin{align*}
		A_3(n)= \dfrac{ \alpha_9\sin (2n+1) \frac{9\pi}{20}}{\sin \frac{9\pi}{20}}- \dfrac{\alpha_1\sin (2n+1) \frac{\pi}{20}}{\sin \frac{\pi}{20}}.
	\end{align*}
	Using MAPLE computations, one can see that
	\begin{align*}
		A_3(20m+0) &=\sqrt{10+2\sqrt{5}},  \qquad
		A_3(20m+1)=\sqrt{10+2\sqrt{5}}, \\
		A_3(20m+2)&=\sqrt{10-2\sqrt{5}}+2\sqrt{10+2\sqrt{5}}, \\
		A_3(20m+3)&=\sqrt{10-2\sqrt{5}}+2\sqrt{10+2\sqrt{5}}, \\
		A_3(20m+4)&=2\sqrt{10-2\sqrt{5}}+2\sqrt{10+2\sqrt{5}}, \\
		A_3(20m+5)&=2\sqrt{10-2\sqrt{5}}+2\sqrt{10+2\sqrt{5}}, \\
		A_3(20m+6)&=\sqrt{10-2\sqrt{5}}+2\sqrt{10+2\sqrt{5}}, \\
		A_3(20m+7)&=\sqrt{10-2\sqrt{5}}+2\sqrt{10+2\sqrt{5}}, \\
		A_3(20m+8)&=\sqrt{10+2\sqrt{5}}, \qquad
		A_3(20m+9)=\sqrt{10+2\sqrt{5}}, \\
		A_3(20m+10)&=-\sqrt{10+2\sqrt{5}},\quad
		A_3(20m+11)=-\sqrt{10+2\sqrt{5}},\\
		A_3(20m+12)&=-\left(\sqrt{10-2\sqrt{5}}+2\sqrt{10+2\sqrt{5}}\right),\\
		A_3(20m+13)&=-\left(\sqrt{10-2\sqrt{5}}+2\sqrt{10+2\sqrt{5}}\right) ,\\
		A_3(20m+14)&=-\left(2\sqrt{10-2\sqrt{5}}+2\sqrt{10+2\sqrt{5}}\right),\\
		A_3(20m+15)&=-\left(2\sqrt{10-2\sqrt{5}}+2\sqrt{10+2\sqrt{5}}\right), \\
		A_3(20m+16)&=-\left(\sqrt{10-2\sqrt{5}}+2\sqrt{10+2\sqrt{5}}\right),\\
		A_3(20m+17)&=-\left(\sqrt{10-2\sqrt{5}}+2\sqrt{10+2\sqrt{5}}\right), \\
		A_3(20m+18)&= -\left(\sqrt{10+2\sqrt{5}}\right),\quad
		A_3(20m+19) =-\left(\sqrt{10+2\sqrt{5}}\right).
	\end{align*}
	Thus,
	\begin{align*}
		&\sum_{n=0}^{\infty} (-1)^n A_3(n) q^{\frac{(2n+1)^2}{8}} = \left(\sqrt{10+2\sqrt{5}}\right)\left\{\sum_{m=0}^{\infty} q^{\frac{(40m+1)^2}{8}}-  \sum_{m=0}^{\infty} q^{\frac{(40m+3)^2}{8}}\right\} \\
		&+\left(\sqrt{10-2\sqrt{5}}+2\sqrt{10+2\sqrt{5}}\right) \left\{\sum_{m=0}^{\infty} q^{\frac{(40m+5)^2}{8}}  - \sum_{m=0}^{\infty} q^{\frac{(40m+7)^2}{8}} \right\}\\
		&+\left(2\sqrt{10-2\sqrt{5}}+2\sqrt{10+2\sqrt{5}}\right) \left\{\sum_{m=0}^{\infty} q^{\frac{(40m+9)^2}{8}}
		-\sum_{m=0}^{\infty}q^{\frac{(40m+11)^2}{8}}\right\}
			\end{align*}
		\begin{align*} 
		&+\left(\sqrt{10-2\sqrt{5}}+2\sqrt{10+2\sqrt{5}}\right) \left\{\sum_{m=0}^{\infty} q^{\frac{(40m+13)^2}{8}}-\sum_{m=0}^{\infty}  q^{\frac{(40m+15)^2}{8}}\right\}\\
		&+\left(\sqrt{10+2\sqrt{5}}\right) \left\{\sum_{m=0}^{\infty}  q^{\frac{(40m+17)^2}{8}} -\sum_{m=0}^{\infty}  q^{\frac{(40m+19)^2}{8}}-\sum_{m=0}^{\infty}  q^{\frac{(40m+21)^2}{8}}+\sum_{m=0}^{\infty}  q^{\frac{(40m+23)^2}{8}}\right\}\\
		&-\left(\sqrt{10-2\sqrt{5}}+2\sqrt{10+2\sqrt{5}}\right)\left\{\sum_{m=0}^{\infty}  q^{\frac{(40m+25)^2}{8}} 
		- \sum_{m=0}^{\infty}  q^{\frac{(40m+27)^2}{8}}\right\} \\
		&-\left(2\sqrt{10-2\sqrt{5}}+2\sqrt{10+2\sqrt{5}}\right) \left\{\sum_{m=0}^{\infty}  q^{\frac{(40m+29)^2}{8}}
		- \sum_{m=0}^{\infty}  q^{\frac{(40m+31)^2}{8}}\right\}\\
		&-\left(\sqrt{10-2\sqrt{5}}+2\sqrt{10+2\sqrt{5}}\right)\left\{ \sum_{m=0}^{\infty}   q^{\frac{(40m+33)^2}{8}}
		-\sum_{m=0}^{\infty}   q^{\frac{(40m+35)^2}{8}}\right\}\\
		&-\left(\sqrt{10+2\sqrt{5}}\right)\left\{ \sum_{m=0}^{\infty}   q^{\frac{(40m+37)^2}{8}}-\sum_{m=0}^{\infty}   q^{\frac{(40m+39)^2}{8}}\right\}.
	\end{align*}
	Changing $m$ to $-m$ and $m$ to $m-1$ in the first and last ten summations respectively, above equation deduces to 
	\begin{align*}
		&\sum_{n=0}^{\infty} (-1)^n A_3(n) q^{\frac{(2n+1)^2}{8}} = \left(\sqrt{10+2\sqrt{5}}\right)\left\{\sum_{m=-\infty}^{\infty} q^{\frac{(40m-1)^2}{8}}-  \sum_{m=-\infty}^{\infty} q^{\frac{(40m-3)^2}{8}} \right. \\ & \left. +\sum_{m=-\infty}^{\infty}  q^{\frac{(40m-17)^2}{8}} -\sum_{m=-\infty}^{\infty}  q^{\frac{(40m-19)^2}{8}}\right\} 
		+\left(\sqrt{10-2\sqrt{5}}+2\sqrt{10+2\sqrt{5}}\right)\\ & \times \left\{\sum_{m=-\infty}^{\infty} q^{\frac{(40m-5)^2}{8}}  - \sum_{m=-\infty}^{\infty} q^{\frac{(40m-7)^2}{8}}+\sum_{m=-\infty}^{\infty} q^{\frac{(40m-13)^2}{8}}-\sum_{m=-\infty}^{\infty}  q^{\frac{(40m-15)^2}{8}} \right\}\\
		&+\left(2\sqrt{10-2\sqrt{5}}+2\sqrt{10+2\sqrt{5}}\right) \left\{\sum_{m=-\infty}^{\infty} q^{\frac{(40m-9)^2}{8}}
		-\sum_{m=-\infty}^{\infty}q^{\frac{(40m-11)^2}{8}}\right\}.
	\end{align*}
	Using the definition of $f(a,b)$, the above expression reduces to 
	\begin{align}\label{3.1.11}
		\nonumber	&\sum_{n=0}^{\infty} (-1)^n A_3(n) q^{\frac{(2n+1)^2}{8}} = \left(\sqrt{10+2\sqrt{5}}\right)q^{1/8}\left\{f(q^{190},q^{210})-  q f(q^{170}, q^{230})   \right. \\ & \left.  \nonumber +q^{36} f(q^{30}, q^{370})-q^{45} f(q^{10},q^{390})\right\} 
		+\left(\sqrt{10-2\sqrt{5}}+2\sqrt{10+2\sqrt{5}}\right)q^{1/8}\\& \nonumber \times  \left\{q^3 f(q^{150}, q^{250})  - q^{6} f(q^{130}, q^{270})+q^{21} f(q^{70}, q^{330}) -q^{28}f(q^{50},q^{350})\right\}\\
		&+\left(2\sqrt{10-2\sqrt{5}}+2\sqrt{10+2\sqrt{5}}\right) q^{1/8}\left\{q^{10} f(q^{110}, q^{290})
		-q^{15} f(q^{90},q^{310})\right\}.
	\end{align}
	Using \eqref{fab1}, \eqref{fab2}, \eqref{fab3}, \eqref{fab4} and \eqref{fab5} in \eqref{3.1.11}, 
	\begin{align}\label{3.1.12}
		\nonumber	&\sum_{n=0}^{\infty} (-1)^n A_3(n) q^{\frac{(2n+1)^2}{8}} = \left(\sqrt{10+2\sqrt{5}}\right)q^{1/8}\left\{f(-q^{45},-q^{55})-qf(-q^{35},-q^{65}) \right\} \\ & \nonumber
		+\left(\sqrt{10-2\sqrt{5}}+2\sqrt{10+2\sqrt{5}}\right)q^{1/8}  \left\{q^3f(-q^{25}, -q^{75})-q^6f(-q^{15},-q^{85})\right\}\\
		& +\left(2\sqrt{10-2\sqrt{5}}+2\sqrt{10+2\sqrt{5}}\right) q^{81/8} f(-q^5,-q^{95}).
	\end{align}
	Utilizing \eqref{3.1.12} in \eqref{3.1.10}, one can see that
	\begin{align}\label{3.1.13}
		\nonumber &	\alpha_9	\Omega_9(q) - \alpha_1 \Omega_1(q)
		= \dfrac{q^{-1/12}}{\eta(\tau)}  \left[\left(\sqrt{10+2\sqrt{5}}\right)q^{1/8}\left\{f(-q^{45},-q^{55})-qf(-q^{35},-q^{65}) \right\} \right. \\ & \left. \nonumber
		+\left(\sqrt{10-2\sqrt{5}}+2\sqrt{10+2\sqrt{5}}\right)q^{1/8}  \left\{q^3f(-q^{25}, -q^{75})-q^6f(-q^{15},-q^{85})\right\} \right. \\ & \left.+\left(2\sqrt{10-2\sqrt{5}}+2\sqrt{10+2\sqrt{5}}\right) q^{81/8} f(-q^5,-q^{95})\right].
	\end{align}
	Multiplying on both the sides of \eqref{3.1.13} by \eqref{prodK},
	\begin{align}\label{3.1.14}
		\nonumber 	\alpha_9\Omega_9(q)\prod_{k=1}^{9} \Omega_k &- \alpha_1  \Omega_1(q)\prod_{k=1}^{9} \Omega_k
		= \dfrac{q^{-5/8}\eta(20\tau)}{\eta(\tau)\eta(2\tau)}\\ & \nonumber \times \left[\left(\sqrt{10+2\sqrt{5}}\right) \left\{f(-q^{45},-q^{55})-qf(-q^{35},-q^{65}) \right\} \right. \\ & \left. \nonumber
		+\left(\sqrt{10-2\sqrt{5}}+2\sqrt{10+2\sqrt{5}}\right)   \left\{q^3f(-q^{25}, -q^{75})-q^6f(-q^{15},-q^{85})\right\} \right. \\ & \left.+\left(2\sqrt{10-2\sqrt{5}}+2\sqrt{10+2\sqrt{5}}\right) q^{10} f(-q^5,-q^{95})\right].
	\end{align}
	Noting the fact that $\alpha_5=0$ and $(-q^2;q^2)_\infty = q^{-2/24}\dfrac{\eta(4\tau)}{\eta(2\tau)}$, upon simplification and changing $q$ to $q^{1/5}$ throughout \eqref{3.1.14}, the identity \eqref{O99-O11} is obtained.  The proofs of the identities \eqref{O99+O11}, \eqref{O77+O33} and \eqref{O77-O33} are similar to the above and hence are omitted. 				
\end{proof}

\section{Lambert series identities related to $T_1(q)$ and $T_2(q)$}
The theory of basic hypergeometric series facilitates a natural framework for studying Ramanujan's continued fractions and related $q$-series identities. Classical summation and transformation formulas reveal deep connections between continued fractions, theta functions, Lambert series, and modular forms. Ramanujan's bilateral summation formula $_1 \psi _1$ acquires a prominent role because of its application in deriving Lambert series expansions and theta-function identities.
Ramanujan's continued fractions of various orders exhibit rich arithmetic and analytic structures, many of which can be understood through hypergeometric techniques. In this section, we explore new connections between the continued fraction of order twenty and the Lambert series. 
\begin{theorem}
	For $|q|<1$,  the following identities hold:
	\begin{align}\label{E1}
		\sum_{\substack{n=1 \\ n\equiv 1 \pmod 2}}^{\infty} \dfrac{q^{3n} +q^{7n}}{1-q^{20n}} - \sum_{\substack{n=1 \\ n\equiv 1\pmod 2}}^{\infty} \dfrac{q^{13n} +q^{17n}}{1-q^{20n}} = \dfrac{\eta^4(40\tau)}{\eta^2(20\tau)} \left[\frac{1}{T_1(q^2)}+ T_1(q^2)\right]
	\end{align}and
	\begin{align}\label{E2}
		\sum_{\substack{n=1 \\ n\equiv 1 \pmod 2}}^{\infty} \dfrac{q^{n} +q^{9n}}{1-q^{20n}} - \sum_{\substack{n=1 \\ n\equiv 1 \pmod 2}}^{\infty} \dfrac{q^{11n} +q^{19n}}{1-q^{20n}} = \dfrac{\eta^4(40\tau)}{\eta^2(20\tau)} \left[\frac{1}{T_2(q^2)}+ T_2(q^2)\right].
	\end{align}
\end{theorem}
\begin{proof}
	Note that
	\begin{align*}
		\sum_{\substack{n=1 \\ n\equiv 1 \pmod 2}}^{\infty} \dfrac{q^{3n} +q^{7n}}{1-q^{20n}} &- \sum_{\substack{n=-\infty \\ n\equiv 1 \pmod 2}}^{-1} \dfrac{q^{-13n} +q^{-17n}}{1-q^{-20n}} \\
		&= \sum_{n=-\infty}^{\infty} \dfrac{q^{6n+3}}{1-q^{40n+20}} +\sum_{n=-\infty}^{\infty} \dfrac{q^{14n+7}}{1-q^{40n+20}}.
	\end{align*}
	Utilizing Ramanujan's $_1 \psi _1$ summation formula \cite[Entry 17, p. 32]{Adiga_1985},
	\begin{align}\label{e1}
		\sum_{n=-\infty}^{\infty} \dfrac{z^n}{1-xq^n} = \dfrac{(xz,q/xz,q,q;q)_\infty}{(x,q/x,z,q/z;q)_\infty}, \qquad |q|<|z|<1,
	\end{align}
	the following is deduced.
	\begin{align}\label{e2}
		\nonumber	\sum_{\substack{n=1 \\ n\equiv 1 \pmod 2}}^{\infty} &\dfrac{q^{3n} +q^{7n}}{1-q^{20n}} - \sum_{\substack{n=-\infty \\ n\equiv 1 \pmod 2}}^{-1} \dfrac{q^{-13n} +q^{-17n}}{1-q^{-20n}} \\&= \dfrac{(q^{40};q^{40})_\infty^2}{(q^{20};q^{40})_\infty^2} \left\{q^3 \dfrac{(q^{14},q^{26};q^{40})_\infty}{(q^6,q^{34};q^{40})_\infty} + q^7 \dfrac{(q^{6},q^{34};q^{40})_\infty}{(q^{14},q^{26};q^{40})_\infty}\right\}.
	\end{align}
	Using the definition of $\eta(\tau)$ and \eqref{T1} and \eqref{T2} in \eqref{e2}, identity \eqref{E1} is obtained. The proof of the identity \eqref{E2} is similar to that of \eqref{E1} and hence is omitted. 
\end{proof}
The left-hand sides of the above identities can be viewed as the generalized Lambert series associated with the residue classes modulo 20.
Let 
\begin{align*}
	F_1(q):= \sum_{\substack{n\geq1 \\ n\equiv 1 \pmod 2}} \dfrac{q^{3n}+q^{7n}-q^{13n}-q^{17n}}{1-q^{20n}}
\end{align*} 
$F_1(q)$ admits the Fourier series expansion:
\begin{align*}
	F_1(q)= \sum_{N\geq 1} a(N) q^N
\end{align*} 
where 
\begin{align*}
	a(N) = \sum_{\substack{d|N \\ d \,\, odd}} \chi(d),
\end{align*}
with 
\begin{align*}
	\chi(d)=
	\begin{cases}
		1, & d\equiv 3,7 \pmod{20},\\
		-1, & d\equiv 13,17 \pmod{20},\\
		0, & \text{otherwise}.
	\end{cases}
\end{align*}
Note that 
\begin{align*}
	q \dfrac{d }{dq} F_1(q) = \sum_{N\geq 1} N a(N) q^N.
\end{align*} is a divisor weighted Lambert series of Eisenstein type.
Similar argument holds for   
\begin{align*}
	F_2(q):= \sum_{\substack{n\geq1 \\ n\equiv 1 \pmod 2}} \dfrac{q^{n}+q^{9n}-q^{11n}-q^{19n}}{1-q^{20n}}.
\end{align*}
Thus, the logarithmic derivatives of $F_1(q)$ and $F_2(q)$ yield  divisor-weighted arithmetic functions having coefficients analogous to Eisenstein coefficients.

\section{Lambert series identities modulo twenty}
In this section, using Jacobi's theta function $\theta_1$ four Lambert series identities associated with the residue classes modulo twenty are established. Differentiating either sides of \eqref{t1}, and substituting $z=0$, we obtain
\begin{align*}
	\theta_1^{'}(0|\tau) =2q^{1/8} (q;q)_\infty^3,
\end{align*} 
where $\theta_1^{'}$ represents the partial derivative of $\theta_1$ with respect to $z$. The following lemma is useful in proving the main results of this section.
\begin{lemma}
	We have 
	\begin{align}\label{Es1}
		\nonumber	\sum_{n=1}^{\infty} &\dfrac{q^n-q^{9n}-q^{11n}+q^{19n}}{1-q^{20n}} \sin2nz \\
		&= -\dfrac{\theta_1(0|20\tau) \theta_1(2z|20\tau) \theta_1(10\pi \tau |20\tau) \theta_1(8\pi \tau | 20\tau)}{4 \theta_1(z-\pi \tau | 20 \tau) \theta_1(z+\pi \tau | 20 \tau) \theta_1(z-9\pi \tau | 20 \tau) \theta_1(z+9\pi \tau | 20 \tau)}
	\end{align} and 
	\begin{align}\label{Es2}
		\nonumber	\sum_{n=1}^{\infty} &\dfrac{q^{3n}-q^{7n}-q^{13n}+q^{17n}}{1-q^{20n}} \sin2nz \\
		&= -\dfrac{\theta_1(0|20\tau) \theta_1(2z|20\tau) \theta_1(10\pi \tau |20\tau) \theta_1(4\pi \tau | 20\tau)}{4 \theta_1(z-3\pi \tau | 20 \tau) \theta_1(z+3\pi \tau | 20 \tau) \theta_1(z-7\pi \tau | 20 \tau) \theta_1(z+7\pi \tau | 20 \tau)}.
	\end{align}
\end{lemma}
\begin{proof}
	For convenience, let $J(z|\tau)$ represent the logarithmic derivative of $\theta_1$ with respect to $z$. Differentiating \eqref{t1} logarithmically, with respect to $z$, upon simplification one can obtain the following identity.
	\begin{align*}
		J\left(z+\frac{\pi \tau}{2}|\tau\right) = -i+4 \sum_{n=1}^{\infty} \dfrac{q^{n/2}}{1-q^n} \sin 2nz.
	\end{align*}
	Replacing $\tau$ by $20\tau$,  the above expression deduces to
	\begin{align*}
		J\left(z+10\pi \tau |20\tau\right) = -i+4 \sum_{n=1}^{\infty} \dfrac{q^{10n}}{1-q^{20n}} \sin 2nz.
	\end{align*}
	Changing $z$ to $z-9\pi \tau$, the above expression deduces to
	\begin{align*}
		J\left(z+\pi \tau|20\tau\right) = -i+4 \sum_{n=1}^{\infty} \dfrac{q^{10n}}{1-q^{20n}} \sin 2n(z-9\pi \tau).
	\end{align*}
	Again, replacing $z$ by $-z$, the above equation becomes
	\begin{align*}
		J\left(z-\pi \tau|20\tau\right) = i+4 \sum_{n=1}^{\infty} \dfrac{q^{10n}}{1-q^{20n}} \sin 2n(z+9\pi \tau).
	\end{align*}
	Summing up the previous two expressions and utilizing the following trigonometric identity
	\begin{align*}
		\sin 2n(z+9\pi\tau) +\sin 2n(z-9\pi\tau) = 2 \cos 18n\pi \tau \sin 2nz,
	\end{align*} in the resulting equation, the following expression is obtained:
	\begin{align*}
		J(z+\pi\tau|20\tau)+	J(z-\pi\tau|20\tau)= 4\sum_{n=1}^{\infty} \dfrac{q^n+q^{9n}}{1-q^{20n}} \sin2nz.
	\end{align*}
	Similarly, one can see that
	\begin{align*}
		J(z+9\pi\tau|20\tau)+	J(z-9\pi\tau|20\tau)= 4\sum_{n=1}^{\infty} \dfrac{q^{9n}+q^{11n}}{1-q^{20n}} \sin2nz.
	\end{align*}
	Combining the above two equations, the following is evident.
	\begin{align}\label{5.1.1}
		\nonumber	&4\sum_{n=1}^{\infty}\dfrac{q^n-q^{9n}-q^{11n}+q^{19n}}{1-q^{20n}} \sin 2nz \\
		&= J(z+\pi\tau|20\tau)+	J(z-\pi\tau|20\tau)-J(z+9\pi\tau|20\tau)-	J(z-9\pi\tau|20\tau).
	\end{align}
	Note that 
	\begin{align*}
		J(y_1|\tau) &+J(y_2|\tau)+	J(y_3|\tau)- J(y_1+y_2+y_3|\tau) \\
		&= \dfrac{\theta_1^{'}(0|\tau) \theta_1(y_1+y_2|\tau) \theta_1(y_2+y_3|\tau) \theta_1(y_1+y_3|\tau)  }{\theta_1(y_1|\tau)  \theta_1(y_2|\tau) \theta_1(y_3|\tau)  \theta_1(y_1+y_2+y_3|\tau)}.
	\end{align*}
	Replacing $\tau$ by $20\tau$ and substituting $y_1=z-\pi\tau$, $y_2=z+\pi\tau$, and $y_3=-z+9\pi\tau$, the above expression is deduced to
	\begin{align}\label{5.1.2}
		\nonumber	J(z+\pi\tau|20\tau) &+	J(z-\pi\tau|20\tau)-J(z+9\pi\tau|20\tau)-	J(z-9\pi\tau|20\tau) \\
		&= \dfrac{\theta_1^{'}(0|20\tau) \theta_1(2z|20\tau) \theta_1(10\pi\tau|20\tau) \theta_1(8\pi\tau|20\tau)  }{\theta_1(z-\pi\tau|20\tau)  \theta_1(z+\pi\tau|20\tau) \theta_1(z-9\pi\tau|20\tau)  \theta_1(z+9\pi\tau|20\tau)}.
	\end{align}
	Combining \eqref{5.1.1} and \eqref{5.1.2}, identity \eqref{Es1} is established. The proof of the identity \eqref{Es2} is similar and hence is omitted.
\end{proof}

\begin{theorem}
	For $|q|<1$, the following identities hold:
	\begin{align}\label{Es3}
		\sum_{n=1}^{\infty} \dfrac{n(q^n-q^{9n}-q^{11n}+q^{19n})}{1-q^{20n}} 
		= \dfrac{q(q^{20};q^{20})_\infty^2 (q^{10};q^{10})^2_\infty (q^8,q^{12};q^{20})_\infty}{(q,q^9,q^{11},q^{19};q^{20})^2_\infty}
	\end{align} and 
	\begin{align}\label{Es4}
		\sum_{n=1}^{\infty} \dfrac{n(q^{3n}-q^{7n}-q^{13n}+q^{17n})}{1-q^{20n}} 
		= \dfrac{q^3(q^{20};q^{20})_\infty^2 (q^{10};q^{10})^2_\infty (q^4,q^{16};q^{20})_\infty}{(q^3,q^7,q^{13},q^{17};q^{20})^2_\infty}.
	\end{align}
\end{theorem}
\begin{proof}
	Dividing both the sides of the equation \eqref{Es1} by $z$ and letting $z \rightarrow 0$, one can see that
	\begin{align*}
		\sum_{n=1}^{\infty} \dfrac{n(q^n-q^{9n}-q^{11n}+q^{19n})}{1-q^{20n}}  = -\dfrac{\theta_1^{'}(0|20\tau)^2 \theta_1(8\pi\tau|20\tau) \theta_1(10\pi \tau |20\tau) }{4 \theta_1^2(\pi \tau | 20 \tau)  \theta_1^2(9\pi \tau | 20 \tau)}.
	\end{align*}
	Using the definition of $\theta_1$ in the right-hand side of the above equation, \eqref{Es3} is established.  The proof of the equation \eqref{Es4} is similar and is omitted.	
\end{proof}

\begin{theorem}
	For $|q|<1$, the following identities hold:
	\begin{align}\label{Es5}
		\nonumber 	\sum_{n=1}^{\infty} \left(\dfrac{n}{3}\right) &\dfrac{(q^n-q^{9n}-q^{11n}+q^{19n})}{1-q^{20n}} 
		\\ &	= \dfrac{q (q^{10};q^{10})^2_\infty (q^{60};q^{60})_\infty (q,q^8,q^9,q^{11},q^{12},q^{19};q^{20})_\infty}{(q^{20};q^{20})_\infty(q^3,q^{27},q^{33},q^{57};q^{60})_\infty}
	\end{align} and 
	\begin{align}\label{Es6}
		\nonumber 	\sum_{n=1}^{\infty} \left(\dfrac{n}{3}\right) &\dfrac{(q^{3n}-q^{7n}-q^{13n}+q^{17n})}{1-q^{20n}} 
		\\ &	= \dfrac{q^3 (q^{10};q^{10})^2_\infty (q^{60};q^{60})_\infty (q^3,q^4,q^7,q^{13},q^{16},q^{17};q^{20})_\infty}{(q^{20};q^{20})_\infty(q^9,q^{21},q^{39},q^{51};q^{60})_\infty}.
	\end{align}
\end{theorem}

\begin{proof}
	Note that 
	\begin{align*}
		\sin \dfrac{2n\pi}{3} = \dfrac{\sqrt{3}}{2} \left(\dfrac{n}{3}\right), \qquad \theta_1\left(\dfrac{2\pi}{3}|\tau\right) = \sqrt{3} q^{1/8} (q^3;q^3)_\infty,
	\end{align*} where $\left(\dfrac{.}{p}\right)$ denotes the Legendre symbol modulo $p$.\\
	Setting $z=\dfrac{\pi}{3}$, the equation \eqref{Es1} deduces to
	\begin{align}\label{5.2.1}
		\nonumber&\sum_{n=1}^{\infty} \left(\dfrac{n}{3}\right)\dfrac{(q^n-q^{9n}-q^{11n}+q^{19n})}{1-q^{20n}} 
		\\ &=- \dfrac{q^{5/2} \theta_1^{'} (0|20\tau) \theta_1(10\pi\tau|20\tau) \theta_1(8\pi\tau|20\tau) (q^{60};q^{60})_\infty}{2\theta_1(\pi/3+\pi\tau|20\tau) \theta_1(\pi/3-\pi\tau|20\tau) \theta_1(\pi/3+9\pi\tau|20\tau) \theta_1(\pi/3-9\pi\tau|20\tau)}.
	\end{align}
	From \cite[eq.(3.1)]{Liu2005},
	\begin{align*}
		\theta_1\left(\frac{\pi}{3}-z|\tau\right) \theta_1\left(\frac{\pi}{3}+z|\tau\right) = \dfrac{(q;q)_\infty^3 \theta_1(3z|3\tau)}{(q^3;q^3)_\infty \theta_1(z|\tau)}.
	\end{align*}
	Using the above equation in \eqref{5.2.1}, \eqref{Es5} is obtained. The proof of \eqref{Es6} is similar and hence is omitted.
\end{proof}

\section{Conclusion}
We have established several new identities connecting Ramanujan's continued fractions of orders five, ten, and twenty, using the product expansion of theta function $\theta_1$. The methods developed in this work suggest a systematic framework for deriving identities connecting continued fractions of various orders through theta function evaluations. We have derived Lambert series identities associated with continued fractions of order twenty. Further, we have also established Lambert series identities associated to residue classes modulo twenty using theta functions.	The identities established here enrich the existing theory of Ramanujan's continued fractions, providing deep insight into the underlying modular structure. The methods described in this article demonstrate the systematic utilization of theta function techniques to generate Lambert series identities. Thus, results presented here contribute to the broader areas of $q$-series, modular form, and special functions inspired by Ramanujan's work.  

\subsection*{Acknowledgments}
The authors acknowledge the guidance of Dr. Vasuki K. R., Department of Studies in Mathematics, University of Mysore, Manasagangotri, Mysore-570 006, India in improving the quality of this research work. The corresponding author gratefully acknowledges the support of DST-INSPIRE, Department of Science and Technology, Government of India, India for INSPIRE fellowship [DST/INSPIRE/03/2022/004970].

		\addcontentsline{toc}{section}{References}
	\bibliographystyle{abbrv}

\end{document}